\documentclass{amsart}

\usepackage{amsmath,amsfonts,amssymb,amsthm,amscd,latexsym,euscript,bbm}

   \newtheorem{lemma}{Lemma}[section]  
   \newtheorem{satz}[lemma]{Theorem}         

\theoremstyle{definition}
\newtheorem{Def}[lemma]{Definition}                                          

    \newtheorem{Bem}[lemma]{Remark}        

\usepackage[all]{xy}

\newenvironment{beweis}{\noindent\textbf{Proof:}}{\par \hfill $\Box$ \hspace{1cm} \par}

\newenvironment{proofof}[1]{\renewcommand\proofname{#1} {\noindent\bf{Proof of #1:}}}%
{\par  \hfill {\bf{$\Box $}} \hspace{1cm}       \par}

\newcounter{zahl}%
\newenvironment{punkt}{\begin{list}{{\rm{(\roman{zahl})}}}%
    {\usecounter{zahl}%
     \setlength{\leftmargin}{0pt} \setlength{\itemindent}{4pt} \setlength{\topsep}{2pt} \setlength{\parsep}{2pt} }}%
    {\end{list}}%

\newcommand{\co}{\colon\thinspace}
\newcommand{\mc}[1]{\ensuremath{\mathcal{#1}}}
\newcommand{\aus}{\raisebox{1pt}{\ensuremath{\,{\scriptstyle\in}\,}}}%

\newcommand{\toh}[1]{\ensuremath{\stackrel{#1}{\rightarrow}}}%
\newcommand{\colim}{\operatorname*{colim}}
\newcommand{\sk}{\operatorname{sk}}
\newcommand{\cosk}{\operatorname{cosk}}
\newcommand{\frei}{\,\_\!\_\,}
\newcommand{\id}{\operatorname{id}}
\newcommand{\Tot}{\operatorname{Tot}}
\newcommand{\inj}{\ensuremath{\hookrightarrow}}

\newcommand{\hel}{\ensuremath{\raisebox{1pt}{[}}}
\newcommand{\her}{\ensuremath{\raisebox{1pt}{]}}}
\newcommand{\hrl}{\ensuremath{\raisebox{1pt}{(}}}
\newcommand{\hrr}{\ensuremath{\raisebox{1pt}{)}}}
\newcommand{\ho}[1]{\ensuremath{H\!o({#1})}}

\newcommand{\bth}{\raisebox{1pt}{\ensuremath{\,\scriptstyle \geq\,}}}%
\newcommand{\sth}{\raisebox{1pt}{\ensuremath{\,\scriptstyle \leq\,}}}%
\newcommand{\uber}{\raisebox{1pt}{\ensuremath{\,{\scriptstyle > }\,}}}%
\newcommand{\ent}{\raisebox{1pt}{\ensuremath{\,{\scriptstyle \subseteq }\,}}}%
\newcommand{\pullback}{\displaystyle \cdot\!\!\lrcorner}%
\newcommand{\winkel}{\ulcorner\hspace{-3pt}\raisebox{1pt}{\ensuremath{\cdot}}}%
\newcommand{\pushout}{\raisebox{-3pt}{\ensuremath{\displaystyle \winkel}}}%

\newcommand{\wt}[1]{\widetilde{#1}}
\newcommand{\stabhom}[2]{\ensuremath{\hel#1,#2\her}}%
\newcommand{\Hom}[3]{\ensuremath{{\rm Hom}_{#1}\hrl#2,#3\hrr}}%
\newcommand{\naturalpi}[3]{\ensuremath{\pi_{#1}^{\natural}\hrl #2,#3\hrr}}%

\newcommand{\map}{\ensuremath{{\rm map}}}
\renewcommand{\hom}{\ensuremath{{\rm hom}}}

\newcommand{\diagr}[1]{ \begin{equation*} \xymatrix{#1} \end{equation*}}%

\numberwithin{equation}{section}        
\newcommand{\Ref}[1]{\hrl\ref{#1}\hrr}        

\begin{document}

\SelectTips{cm}{10}

\title  [Truncated resolution structures]
        {Truncated resolution model structures}
\author{Georg Biedermann}

\address{Department of Mathematics, Middlesex College, The University of Western Ontario, London, Ontario N6A 5B7, Canada}

\email{gbiederm@uwo.ca}



\date{\today}
\dedicatory{}
\commby{}

\begin{abstract}
Using the dual of Bousfield-Friedlander localization we colocalize resolution model structures on cosimplicial objects over a left proper model category to get truncated resolution model structures. These are useful to study realization and moduli problems in algebraic topology.
\end{abstract}

\maketitle

\section{Introduction} 
\label{section:intro}

Resolution model structures were first introduced in \cite{DKSt:E2} and later studied in \cite{DKSt:bigraded} and \cite{BlDG:pi-algebra} to attack the realization problem for $\Pi$-algebras. Following these tracks a similar resolution model structures was developed in \cite{GoHop:resolution} and used in \cite{GoHop:moduli} to study $A_\infty$- and $E_\infty$-structures on ring spectra. 
In \cite{Bou:cos} a very general and elegant treatment of resolution model structures is given that exhibits the previous ones as special cases. Bousfield in his paper calls these structures the \mc{G}-resolution structures because there is the freedom of choosing an appropriate class of injective models \mc{G} which will be explained in subsection \ref{subsection:G-Struktur}. We will abbreviate this and call the resolution model structures simply \mc{G}-structures.

What these resolution model structures allow us to do is to think of cosimplicial objects as analogues of cochain complexes and of fibrant approximations as analogues of injective resolutions. 
More precisely one defines homotopy and homology of a cosimplicial object with respect to the above mentioned class of injective models \mc{G}. Then weak equivalences are given in terms of either homotopy or homology.

The goal of this work is to describe for each $n\bth 0$ an $n$-truncated resolution model structure. Its weak equivalences are maps that induce isomorphisms of homotopy just up to degree $n$. We will obtain these structures by colocalization with respect to an augmented homotopy-idempotent homotopy functor in the sense of \cite{BF:gamma}. 

We use these truncated structures in \cite{Biedermann:interpol} to study extensively the realization and moduli problem with respect to a nice class of homology theories and to define interpolation categories for them. This article develops the necessary technical framework. It is the first part of my thesis written at the Universit\"at Bonn. I would like to thank my advisor Jens Franke and my coadvisor Stefan Schwede. I am very grateful to Pete Bousfield for several very helpful e-mails as well as for sending me his then unpublished paper \cite{Bou:cos}. I also thank the referee for his comments.

\section{Resolution model structures} 
\label{section:E2-Strukturen}
This section first introduces cosimplicial objects and the Reedy structure and then reviews the theory of resolution model structures.

\subsection{Cosimplicial objects and the Reedy structure}
\label{section:Reedy}
Let \mc{M} be a simplicial model category. Let $c\mc{M}$ be the category of cosimplicial objects over \mc{M}. 
We refer to \cite{GoJar:simp}, \cite{Hir:loc} or \cite{Hov:model} for the necessary background, in particular for the internal simplicial structure, which is compatible with the Reedy structure, and for latching- and matching objects. A quick description of latching objects can be found in \ref{partielle Matching- und Latchingobjekte}. Beware of a degree shift between our matching objects and the ones in \cite{GoJar:simp}.

\begin{Def} \label{Reedy Modellstruktur}
We define the following classes of morphisms that will constitute the {\bf Reedy structure} on $c\mc{M}$. Let $X^\bullet\to Y^\bullet$ be a morphism in $c\mc{M}$. It is called 
\begin{punkt}
    \item a {\bf Reedy equivalence} if for every $s\aus\mathbbm{N}$ the maps $X^s\to Y^s$ are weak equivalences in \mc{M}.
    \item a {\bf Reedy cofibration} if for every $s\aus\mathbbm{N}$ the induced maps 
    $$X^s\sqcup_{L^{s}X^\bullet}L^{s}Y^\bullet\to Y^s$$
are cofibrations in \mc{M}.
    \item a {\bf Reedy fibration} if for every $s\aus\mathbbm{N}$ the induced maps 
    $$X^s\to M^{s}X^\bullet\times_{M^{s}Y^\bullet}Y^s$$
are fibrations in \mc{M}.
\end{punkt}
\end{Def}

The following theorem was proved in \cite{Ree:homotopy}, see also \cite[VII.2.12.]{GoJar:simp} and \cite[15.3.4.]{Hir:loc}. 
\begin{satz}
The category $c\mc{M}$ together with the Reedy structure becomes a model category. It becomes a simplicial model structure if we provide it with the internal simplicial structure.
\end{satz}

\subsection{The external simplicial structure on $c\mc{M}$}
\label{subsection:extern}
The resolution model structures are not compatible with the internal simplicial structure. Here we describe the external simplicial structure, which will be compatible with the resolution structure and its truncated versions.

\begin{Bem} \label{Tensor als Koende}
For $X^\bullet$ in $c\mc{M}$ and $L$ in \mc{S} we can perform the following coend-construction: Let $\bigsqcup_{L_\ell}X^m$ be the coproduct in \mc{M} of copies of $X^m$ indexed by the set $L_\ell$, and view this as a functor $\Delta^{\rm op}\times\Delta\to\mc{M}$. Then we can take the coend
    $$ X^\bullet\otimes_{\Delta}L:=\int^{\Delta}\bigsqcup_{L_\ell}X^m \aus\mc{M} .$$
Explicitly this is given by the coequalizer 
    $$ \bigsqcup_{\ell\to m}\bigsqcup_{L_\ell} X^m \rightrightarrows \bigsqcup_{\ell\bth 0}\bigsqcup_{L_\ell} X^\ell \to  X^\bullet\otimes_{\Delta}L   $$
using the obvious maps induced by $\ell\to m$.
\end{Bem}

We are now ready to describe the functors that will enrich all our model structures to simplicial model categories.

\begin{Def} \label{Externe simpliziale Struktur}
We define a simplicial structure on $c\mc{M}$. Let $K$ be in \mc{S} and $X^\bullet$ and $Y^\bullet$ in $c\mc{M}$, then set
\begin{equation*}
          \hrl X^\bullet\otimes^{\rm ext} K\hrr^n  :=\, X^\bullet\otimes_{\Delta} \hrl K\times\Delta^n\hrr , 
\end{equation*}
where $\times$ denotes the usual product of simplicial sets and $\Delta^n$ is the standard $n$-simplex,
\begin{equation*}
      \hspace*{-2.2cm} \hom^{\rm ext}\hrl K,X^\bullet\hrr^n :=\ \ \prod_{K_n}X^n, 
\end{equation*} 
where the product is taken over the set of $n$-simplices of $K$, and finally
\begin{equation*}
      \hspace*{5mm}    \map^{\rm ext}\hrl X^\bullet,Y^\bullet\hrr_n :=\ \Hom{c\mc{M}}{
X^\bullet\otimes^{\rm ext}\Delta^n}{Y^\bullet}.
\end{equation*} 
We call this the {\bf external (simplicial) structure} on $c\mc{M}$. Note that we do not refer to any simplicial structure of \mc{M}. 
From now on we will usually drop the superscripts.
\end{Def} 

\begin{Bem} 
Using the $\Delta$-coend from remark \ref{Tensor als Koende} we observe that we have isomorphisms:
   $$ X^n \cong  X^\bullet\otimes_{\Delta}\Delta^n \hbox{ and } 
      L^nX^\bullet \cong X^\bullet\otimes_{\Delta}\partial\Delta^n $$
Partial latching objects are obtained by inserting various subcomplexes of $\Delta^n$.
\end{Bem}

\begin{Def} \label{partielle Matching- und Latchingobjekte}
Let $\iota_n$ be the unique non-degenerate $n$-simplex of $\Delta^n$.
For $0\sth k\sth n$ let $\Lambda^n_k\ent\Delta^n$ be the subcomplex spanned by all $d_i\iota_n$ for $i\neq k$. We call this the {\bf\boldmath $k$-horn\unboldmath} on $\Delta^n$.
Let $X^\bullet$ be an object of $c\mc{M}$. Set
    $$ L^n_kX^\bullet:=X^\bullet\otimes_{\Delta}\Lambda^n_k .$$
These objects are called the {\bf partial latching objects} of $X^\bullet$. 
\end{Def}

\begin{Def} \label{äußere Einhängung}\label{aeussere Einhaengung}
For an object $X^\bullet$ in $c\mc{M}$ we define its {\bf\boldmath $s$-th external suspension\unboldmath} $\Sigma_{\rm ext}^sX^\bullet$ by the following pushout diagram:
\diagr{ X^\bullet=X^\bullet\otimes^{\rm ext}\ast \ar[r]\ar[d]\ar@{}[dr]|->>>>{\phantom{xx}\pushout} & X^\bullet\otimes^{\rm ext}\Delta^s/\partial\Delta^s \ar[d]& \\
        \ast \ar[r] &  X^\bullet\wedge^{\rm ext}\Delta^s/\partial\Delta^s\ar@{=}[r] & :\Sigma_{\rm ext}^sX^\bullet }
\end{Def}

\subsection{The \mc{G}-structure on cosimplicial objects}
\label{subsection:G-Struktur}

This subsection recapitulates the results of \cite{Bou:cos}, \cite{DKSt:bigraded} and \cite{DKSt:E2}.
We want to think of cosimplicial objects as resolutions of objects in \mc{M} and we want to identify different resolutions of the same object of \mc{M}. Like in the common situation of complexes over an abelian category we have to construct a certain model structure on $c\mc{M}$ which enables us to compare objects with respect to (co-)homology or (co-)homotopy groups.
So the first thing we have to find is the right notion of homotopy and homology of cosimplicial objects. This was done in \cite{DKSt:bigraded} and \cite{DKSt:E2} working simplicially.
The resulting model structure is called a resolution or $E_2$-model structure.
Here we deal with the dual situation and we distinguish objects by mapping into some class of injective models, instead of e.g. mapping out of spheres in the classical case.
To construct the truncated versions we will only consider the natural homotopy groups from \ref{natural homotopy groups}. Still we put the emphasis on considering both two kinds of homotopy groups, since it is suggested by \cite{DKSt:bigraded} and \cite{GoHop:moduli} and seems natural when considering realization problems.

The following definitions are taken from \cite{Bou:cos} who gave the definitive treatment on resolution model structures. 
\begin{Def} \label{injektive Modelle}
Let \mc{M} be a left proper pointed model category. We call a class \mc{G} of objects in \mc{M} a class of {\bf injective models} if the elements of \mc{G} are fibrant and group objects in the homotopy category \ho{\mc{M}} and if \mc{G} is closed under loops. We reserve the letter \mc{G} for such a class.
\end{Def} 

Now we are going to associate to a cosimplicial object groups that the reader should consider as its cohomology. 
They are contravariant functors on $c\mc{M}$ and depend on two parameters.
In the situation of definition \ref{injektive Modelle} let $X^\bullet$ be an object in $c\mc{M}$ and let $\text{ho}\,\mc{G}$ be the class $\mc{G}$ considered as a full subcategory of $\ho{\mc{M}}=\mc{T}$. 
Let $\stabhom{\frei}{\frei}$ denote the morphisms in \mc{T}. 
Note that $\stabhom{X^\bullet}{G}$ is a simplicial group.
For every $s\bth 0$ we have a functor
\begin{equation}\begin{split}\begin{array}{ccc}
    \text{ho}\,\mc{G} & \to & \text{groups} \\
       G   & \mapsto & \pi_s\stabhom{X^\bullet}{G} .
    \end{array}\end{split} 
\end{equation}
Note that for $s\uber 0$ these groups are actually abelian.
On the other hand we can also consider the pointed simplicial set $\Hom{\mc{M}}{X^\bullet}{G}$, where the constant map $X^0\to G$ of the pointed category \mc{M} serves as basepoint. 
If $X^\bullet$ is Reedy cofibrant then this simplicial set is fibrant. 
It supplies a functor
\begin{equation}\begin{split}\begin{array}{ccc}
    \mc{G} & \to & \text{homotopy group objects} \\
       G   & \mapsto & \Hom{\mc{M}}{X^\bullet}{G} . 
             \end{array}\end{split}
\end{equation}
where $\mc{G}$ is considered as a full subcategory of \mc{M}.
Its homotopy should be thought of as the (co-)homotopy of $X^\bullet$. Observe also the equality:
\begin{equation} \label{Hom=map}
    \Hom{\mc{M}}{X^\bullet}{G}= \map^{\rm ext}\hrl X^\bullet, r^0G\hrr ,
\end{equation}
where $r^0G$ denotes the constant cosimplicial object over $G$.

\begin{Def} \label{natural homotopy groups}
Following \cite{GoHop:moduli} we denote the homotopy groups of these group objects by
    $$ \naturalpi{s}{X^\bullet}{G} := \pi_s\Hom{\mc{M}}{X^\bullet}{G} = \pi_s\map\hrl X^\bullet,r^0G\hrr $$
for $s\bth 0$ and $G\aus\mc{G}$ and call them the {\bf\boldmath natural homotopy groups of $X^\bullet$ with coefficients in $\mc{G}$}. 
Note that $r^0G$ is Reedy fibrant, so, again, these groups have homotopy meaning if $X^\bullet$ is Reedy cofibrant.
\end{Def}

\begin{Bem} 
Obviously the canonical functor $\mc{M}\to\ho{\mc{M}}$ induces a map $\Hom{\mc{M}}{X^\bullet}{G}\to\stabhom{X^\bullet}{G}$ which in turn induces a natural transformation of functors
\begin{equation}\nonumber 
     \naturalpi{s}{X^\bullet}{G}\to \pi_s\stabhom{X^\bullet}{G} .
\end{equation}
This map is called the {\bf Hurewicz map} and was constructed in \cite[7.1]{DKSt:bigraded}. 
One of the main results is \cite[8.1]{DKSt:bigraded} (also \cite[3.8]{GoHop:moduli} and \cite{GoHop:resolution}) that this Hurewicz homomorphism for each $G\aus\mc{G}$ fits into a long exact sequence, the so-called {\bf spiral exact sequence}
\begin{align*}
    ...& \to \naturalpi{s-1}{X^\bullet}{\Omega G}\to\naturalpi{s}{X^\bullet}{G}\to\pi_s\stabhom{X^\bullet}{G}\to \naturalpi{s-2}{X^\bullet}{\Omega G} \to ... \\
    ...& \to \pi_{2}\stabhom{X^\bullet}{G}\to \naturalpi{0}{X^\bullet}{\Omega G}\to \naturalpi{1}{X^\bullet}{G}\to \pi_1\stabhom{X^\bullet}{G} \to 0,
\end{align*}
where $\Omega$ is the loop space functor on \mc{M}, plus an isomorphism 
    $$ \naturalpi{0}{X^\bullet}{G}\cong\pi_0\stabhom{X^\bullet}{G}.  $$
In the construction of the exact sequence we rely on the external simplicial structure, but not on a simplicial structure of \mc{M}. 

As explained in \cite[8.3.]{DKSt:bigraded} or \cite[(3.1)]{GoHop:moduli} these long exact sequences can be spliced together to give an exact couple and an associated spectral sequence
\begin{equation} \label{Spiralspektral}
      \pi_p\stabhom{X^\bullet}{\Omega^q G}\ \Longrightarrow\ \colim_k\naturalpi{k}{X^\bullet}{\Omega^{p+q-k}G} 
\end{equation}
\end{Bem}

\begin{Def} 
We call a map in $c\mc{M}$ a {\bf\boldmath \mc{G}-equivalence} if it induces isomorphisms 
    $$ \pi_s\stabhom{Y^\bullet}{G}\to\pi_s\stabhom{X^\bullet}{G} $$
for all $s\bth 0$ and all $G\aus\mc{G}$.
\end{Def}

\begin{lemma} \label{naturale G-Äquivalenzen}\label{naturale G-Aequivalenzen} 
A map $X^\bullet\to Y^\bullet$ in $c\mc{M}$ is a \mc{G}-equivalence if and only if it induces iso\-mor\-phisms
    $$ \emph{\naturalpi{s}{\wt{Y}^\bullet}{G}}\to\emph{\naturalpi{s}{\wt{X}^\bullet}{G}} $$
for all $s\bth 0$ and all $G\in\mc{G}$ with the canonical basepoint and some Reedy cofibrant approximation $\wt{X}^\bullet\to\wt{Y}^\bullet$.
\end{lemma}

\begin{beweis}
We note first, that it suffices to choose one basepoint since the objects are fibrant homotopy group objects.
Now this follows immediately from the spiral exact sequence by simultaneous induction over the whole class \mc{G} and the five-lemma. Remember that \mc{G} is closed under loops by assumption.
\end{beweis}

To describe the fibrations and cofibrations of the resolution or \mc{G}-model structure on $c\mc{M}$ we have to introduce further definitions.
\begin{Def} \label{G-monisch, G-injektiv}
A map $A\toh{i} B$ in \mc{M} is called {\bf \mc{G}-monic} when $\stabhom{B}{G}\toh{i^*}\stabhom{A}{G}$ is surjective for each $G\aus\mc{G}$. 

An object $I$ is called {\bf\boldmath$\mc{G}$\unboldmath-injective} when $\stabhom{B}{I}\toh{i^*}\stabhom{A}{I}$ is surjective for each \mc{G}-monic map $A\toh{i} B$. 

We call a fibration in \mc{M} a {\bf\boldmath$\mc{G}$\unboldmath-injective fibration} if it has the right lifting property with respect to every \mc{G}-monic cofibration.

We say that \ho{\mc{M}} {\bf has enough \boldmath$\mc{G}$\boldmath-injectives} if each object in \ho{\mc{M}} is the source of a \mc{G}-monic map to a \mc{G}-injective target. We say that \mc{G} is {\bf functorial}, if these maps can be chosen functorially.
\end{Def}

\begin{Def} 
A map $X^\bullet\to Y^\bullet$ will be called a {\bf\boldmath \mc{G}-fibration} if $f:X^n\to Y^n\times_{M^nY^{\bullet}}M^nX^{\bullet}$ is a \mc{G}-injective fibration for $n\bth 0$.
\end{Def}

Note that there is yet another characterization of \mc{G}-equivalences, namely these are the maps $X^\bullet\to Y^\bullet$ that induce weak equivalences
    $$ \map\hrl \wt{Y}^\bullet, r^0G\hrr\to\map\hrl \wt{X}^\bullet, r^0G\hrr, $$
where $\wt{X}^\bullet\to\wt{Y}^\bullet$ is a cofibrant approximation to $X^\bullet\to Y^\bullet$. Using these induced maps we can also define the cofibrations of our resolution model structure.

\begin{Def} 
We call a map $X^{\bullet}\to Y^{\bullet}$ a {\bf\boldmath \mc{G}-cofibration} if it is a Reedy cofibration and the induced map $\map\hrl Y^{\bullet},r^0G\hrr\to\map\hrl X^{\bullet},r^0G\hrr$ is a fibration of simplicial sets for each $G\aus\mc{G}$. 
These three classes of \mc{G}-equivalences, \mc{G}-cofibrations and \mc{G}-fibrations will be called the {\bf\boldmath \mc{G}-structure on $c\mc{M}$\unboldmath}. We denote it by {\boldmath $c\mc{M}^{\mc{G}}$}.
\end{Def}

\begin{lemma} 
For a map $i:A^\bullet\to B^\bullet$ in $c\mc{M}$ the following are equivalent:
\begin{punkt}
    \item The map $i$ is an \mc{G}-cofibration.
    \item The map $i$ is a Reedy cofibration and for every $G\aus\mc{G}$ the induced map
     $$ \emph{\stabhom{B^\bullet}{G}}\to\emph{\stabhom{A^\bullet}{G}}$$
is a fibration of simplicial sets.
    \item The map $i$ is a Reedy cofibration and the induced maps
       $$ A^s\sqcup_{L^{s}_k A^\bullet}L^{s}_k B^\bullet\to B^s$$
are \mc{G}-monic cofibrations for all $s\bth 1$ and $s\bth k\bth 0$.
\end{punkt}
\end{lemma}

\begin{beweis}
The equivalence of (i) and (iii) is proved by adjunction, \cite[3.13.]{Bou:cos}. The equivalence of (ii) and (iii) is \cite[3.15.]{Bou:cos}. 
\end{beweis}

The message that the reader should get from this is that it does not matter whether we use the simplicial abelian groups $\stabhom{X^\bullet}{G}$ or the fibrant $H$-spaces $\map\hrl X^\bullet, r^0G\hrr$ to define the \mc{G}-structure. It is also worth noting, that the \mc{G}-structure just depends on the \mc{G}-injectives.
Here is the main theorem 3.3 from \cite{Bou:cos}.
\begin{satz} \label{G-Struktur}
Let $\mc{M}$ be a left proper pointed model category \mc{M} with a class \mc{G} of injective models and enough \mc{G}-injectives. Then the \mc{G}-structure on the category $c\mc{M}$ becomes a simplicial left proper pointed model category with the external simplicial structure.  
If \mc{G} is functorial and the model structure on \mc{M} has functorial factorizations, then so does $c\mc{M}^{\mc{G}}$.
\end{satz}

\begin{Bem} 
There is another way of describing the natural homotopy groups, namely
    $$ \naturalpi{s}{X^\bullet}{G}=\stabhom{X^\bullet}{\Omega^s_{\rm ext}r^0G}_{\mc{G}}, $$
where $\stabhom{\frei}{\frei}_{\mc{G}}$ denotes the morphisms in $\ho{c\mc{M}^\mc{G}}$ and $\Omega_{\rm ext}$ is the dual construction to \ref{aeussere Einhaengung}.
\end{Bem}

\section{Truncated resolution model structures}
\label{section:n-G-Struktur}

In this subsection we will truncate the \mc{G}-structure. This means that we are going to construct for each $n\bth 0$ a new model structure, called the $n$-\mc{G}-structure, whose weak equivalences are maps that induce isomorphisms of natural homotopy groups just up to degree $n$. We will obtain these structures by colocalization with respect to a coaugmented homotopy-idempotent homotopy functor. 

\subsection{Colocalization \'{a} la Bousfield-Friedlander}
\label{subsection:BF}
An easy method for localizing with respect to an augmented homotopy-idempotent homotopy functor was devised in \cite{BF:gamma} and considerably improved in \cite{Bou:telescopic}. The approach completely dualizes since no small object argument is used.

We would like to point out, that the truncation or colocalization process does not work for the other type of homotopy groups, which the author had to learn painfully while struggling for the right constructions: Never confuse (co-)limits with homotopy (co-)limits!

\begin{Def} \label{Ko-Q-Struktur}
Let \mc{N} be a left proper model category. A {\bf\boldmath co-$Q$-structure} on \mc{N} consists of a functor $Q\co \mc{N}\to\mc{N}$ and a natural transformation $\alpha\co Q\to\id$ satisfying the following axioms:
\begin{punkt}
    \item If $X\to Y$ is a weak equivalence, then so is $QX\to QY$.
    \item For each $X$ in \mc{N}, the maps $\alpha_{QX}$ and $Q\alpha_X\co QQX\to QX$ are weak equivalences.
    \item For a pushout square
\diagr{ A \ar[r]^v\ar[d]_f\ar@{}[dr]|->>{\pushout} & V \ar[d]^g \\ B\ar[r]_w & W }
in \mc{N}, if $f$ is a cofibration between cofibrant objects such that $\alpha_A$, $\alpha_B$ and $Qv$ are weak equivalences, then $Qw$ is also a weak equivalence.
\end{punkt}
\end{Def}

Reminiscent of \cite{BF:gamma} we use the same letter $Q$ and invented the word co-$Q$-structure for the dual concept. The axioms and definitions together with the following theorem simply state, that the colocalization with respect to the class of $Q$-equivalences exist if we can assure left properness of this structure in advance.
\begin{Def} \label{assoz. Ko-Q-Struktur}
Let $X\to Y$ be a map in \mc{N}, then we say, that it is
\begin{punkt}
    \item a {\bf\boldmath $Q$-equivalence} if $QX\to QY$ is a weak equivalence.
    \item a {\bf\boldmath $Q$-fibration} if it is a fibration.
    \item a {\bf\boldmath $Q$-cofibration} if it has the left lifting property with respect to all $Q$-trivial fibrations, i.e. all maps that are fibrations and $Q$-equivalences.
\end{punkt}
The next theorem proves this to be a model structure, which we will call the $Q$-colocal structure and we will write {\boldmath$\mc{N}^Q$} for it.
\end{Def}

\begin{satz} \label{Ko-Q-Modellstruktur}
Let \mc{N} be a left proper model category equipped with a co-$Q$-structure.
The associated model structure $\mc{N}^Q$ is a left proper model structure. A map $i\co A\to B$ in \mc{N} is a $Q$-cofibration if and only if $i$ is a cofibration in \mc{N} and the square
\diagr{ QA \ar[r]^{\alpha_A}\ar[d]_{Qi} & A \ar[d]^{i} \\
        QB \ar[r]_{\alpha_B} & B }
is a homotopy pushout square in \mc{N} with its original structure. If \mc{N} carries additionally a simplicial structure, then the $Q$-colocal structure on \mc{N} is also simplicial.
If \mc{N} is proper or possesses functorial factorization, the same is true for $\mc{N}^Q$.
\end{satz}

The dual form of this theorem is proved as Theorem 9.3. and Theorem 9.7 of \cite{Bou:telescopic} and it builds upon an earlier version in \cite{BF:gamma}. The fact that left properness suffices to prove left properness is noted in \cite[Remark 9.5]{Bou:telescopic}.

\subsection{The $n$-\mc{G}-structure on $c\mc{M}$}
Now we will come to the heart of the matter and truncate the resolution model structures.

\begin{Def} \label{n-G-structure}
For $G\aus\mc{G}$ we denoted by $r^0G$ the constant cosimplicial object over $G$.
Let $Q_n\co c\mc{M}\to c\mc{M}$ be the composition of a Reedy cofibrant replacement functor $\wt{\frei}$ with the functor $\sk_{n+1}$ and let $\alpha_{X^\bullet}$ be the canonical map 
    $$\sk_{n+1}\wt{X}^\bullet\to\wt{X}^\bullet\to X^\bullet.$$
By applying the above constructions we get for each $n\bth 0$ a new model structure on $c\mc{M}$. We call it the {\bf\boldmath $n$-$\mc{G}$-structure\unboldmath} and denote it by {\boldmath $c\mc{M}^{\text{$n$-\mc{G}}}$}. 
Translating definition \ref{assoz. Ko-Q-Struktur} into our special situation we call a map $f\co X^\bullet\to Y^\bullet$
\begin{punkt}
    \item an {\bf\boldmath $n$-$\mc{G}$-equivalence\unboldmath} if the induced map
            $$ \sk_{n+1}\wt{X}^\bullet\to\sk_{n+1}\wt{Y}^\bullet$$
is a \mc{G}-equivalence.
    \item an {\bf\boldmath $n$-\mc{G}-fibration\unboldmath} if it is a \mc{G}-fibration.
    \item an {\bf\boldmath $n$-\mc{G}-cofibration\unboldmath} if it has the left lifting property with respect to all \mc{G}-fibrations that are also $n$-\mc{G}-equivalences.
\end{punkt}
The original \mc{G}-structure can be thought of as the limit for $n=\infty$.
\end{Def}

\begin{satz} \label{abgeschnittene Modellstruktur}
Let $n\bth 0$. For a pointed left proper model category \mc{M} the $n$-\mc{G}-structure on the category $c\mc{M}$ of cosimplicial objects over \mc{M} is a pointed simplicial left proper model structure. If $\mc{M}^\mc{G}$ possesses functorial factorization then so does $c\mc{M}^{n-\mc{G}}$.
\end{satz}

This theorem is a direct consequence of theorem \ref{Ko-Q-Modellstruktur} once we prove that $Q_n=\sk_{n+1}\wt{\frei}$ provides a co-$Q$-structure. 
If $c\mc{M}^\mc{G}$ is right proper, then so is the $n$-\mc{G}-structure, but we do not know good conditions that ensure the right properness of the \mc{G}-structure. However, it happens, see \cite{Biedermann:interpol}.
The rest of this section is devoted to the proof of theorem \ref{abgeschnittene Modellstruktur}. 
A characterization of $n$-\mc{G}-cofibrations is given in \ref{n-G-Kofaserungen}.

\begin{Bem} \label{map und konstante Objekte}
Observe that by \Ref{Hom=map} there is a natural isomorphism
    $$ \map\hrl\sk_{n+1}X^\bullet, r^0G\hrr\cong\cosk_{n+1}\map\hrl X^\bullet, r^0G\hrr .$$
If $X^\bullet$ is Reedy cofibrant then $\map\hrl X^\bullet, r^0G\hrr$ is fibrant. Note also, that for a Kan-complex $W$ the space $\cosk_{n+1}W$ is a model for the $n$-th Postnikov section.
\end{Bem}

\begin{lemma} \label{bis Grad n}
For a map $f\co X^\bullet\to Y^\bullet$ in $c\mc{M}$ the following assertions are equivalent:
\begin{punkt}
    \item $f$ is an $n$-\mc{G}-equivalence.
    \item  For every $G\in\mc{G}$ and all $0\sth s\sth n$ the induced maps
    $$ \emph{\naturalpi{s}{\wt{Y}^\bullet}{G}}\to\emph{\naturalpi{s}{\wt{X}^\bullet}{G}} $$
are isomorphisms, where $\wt{X}^\bullet\to\wt{Y}^\bullet$ is a cofibrant approximation to $X^\bullet\to Y^\bullet$.
\end{punkt}
\end{lemma}

\begin{beweis}
The equivalences follow readily from remark \ref{map und konstante Objekte}.
\end{beweis}

\begin{proofof}{\ref{abgeschnittene Modellstruktur}}
For arbitrary $G\aus\mc{G}$ and $X^\bullet$ in $c\mc{M}$ we compute with \ref{map und konstante Objekte}:
\begin{equation*}
      \naturalpi{s}{Q_nX^\bullet}{G}\cong\pi_s\cosk_{n+1}\map\hrl\wt{X}^\bullet,r^0G\hrr \cong 
                              \left\{ \begin{array}{cl}
                        \naturalpi{s}{\wt{X}^\bullet}{G} &,\hbox{ for } 0\sth s\sth n \\
                                    0                    &,\hbox{ for } s\uber n   
                                      \end{array}
                              \right.
\end{equation*}
Now conditions (i) and (ii) of \ref{Ko-Q-Struktur} are obvious. To prove (iii) we can assume that all objects in the pushout square are Reedy cofibrant by factoring $v$ and $w$ appropriately. We get a pullback diagram of the following form:
\diagr{ \map\hrl W^\bullet, r^0G\hrr \ar[r]\ar[d] & \map\hrl B^\bullet, r^0G\hrr \ar[d] \\
        \map\hrl V^\bullet, r^0G\hrr \ar[r] & \map\hrl A^\bullet, r^0G\hrr \ar@{}[ul]|->>{\pullback}}
We apply $\naturalpi{s}{\frei}{G}$ and using the fact that $\alpha_{A^\bullet}, \alpha_{B^\bullet}$ and $Q_nv$ are equivalences we obtain, that this square is transformed into a pullback square of groups:
    $$ \naturalpi{s}{W^\bullet}{G}\cong\naturalpi{s}{B^\bullet}{G}\times_{\pi_{s}^\natural(A^\bullet,G)}\naturalpi{s}{V^\bullet}{G}\cong 
                              \left\{ \begin{array}{cl}
                        \naturalpi{s}{B^\bullet}{G} &,\hbox{ for } 0\sth s\sth n \\
                                    0                    &,\hbox{ for } s\uber n   
                                      \end{array}
                              \right.$$
\end{proofof}

We will now give some characterizations of the cofibrations of the $n$-\mc{G}-structure to complete the picture. Of course, since we have not changed the class of fibrations, $n$-\mc{G}-trivial cofibrations are the same as \mc{G}-trivial cofibrations.
The characterization of $n$-\mc{G}-cofibrations is analogous to lemma \cite[3.13]{Bou:cos}.

\begin{lemma} \label{n-G-Kofaserungen}
For a map $i\co A^\bullet\to B^\bullet$ in $c\mc{M}$ the following are equivalent:
\begin{punkt}
    \item The map $i$ is an $n$-\mc{G}-cofibration.
    \item The map $i$ is a \mc{G}-cofibration and the induced maps
       $$ \emph{\naturalpi{s}{\wt{B}^\bullet}{G}}\to\emph{\naturalpi{s}{\wt{A}^\bullet}{G}} $$
are isomorphisms for every $G\aus\mc{G}$ and all $s\uber n$, where $\wt{A}^\bullet\to\wt{B}^\bullet$ is a Reedy cofibrant approximation to $i$.
    \item The map $i$ is a \mc{G}-cofibration and for all $s\bth n+2$ the maps
       $$ A^s\sqcup_{L^{s}A^\bullet}L^{s}B^\bullet\to B^s$$
and the map
       $$ L^{n+2}A^\bullet\sqcup_{A^0}B^0\to L^{n+2}B^\bullet$$
are \mc{G}-monic. Here the last map is induced by some map $\ast\to\partial\Delta^{n+2}$.
\end{punkt}
\end{lemma}

\begin{beweis}
Let $i$ be an $n$-\mc{G}-cofibration. Then by definition it has the left lifting property with respect to $n$-\mc{G}-trivial fibrations and in particular it is a Reedy cofibration.
By the characterization of $Q$-cofibrations in theorem \ref{Ko-Q-Modellstruktur} there is the following homotopy pushout square for the map $i$:
\diagr{ \sk_{n+1}\wt{A}^\bullet \ar[r]^-{\alpha_{A^\bullet}}\ar[d]_{\sk_{n+1}\wt{i}}\ar@{}[dr]|->>>{\text{ho-}\pushout} & A^\bullet \ar[d]^i \\
        \sk_{n+1}\wt{B}^\bullet \ar[r]_-{\alpha_{B^\bullet}}    & B^\bullet  }
The same considerations as in the proof of \ref{abgeschnittene Modellstruktur} show that applying the functor $\pi_s\map\hrl\wt{\frei},r^0G\hrr$ yields a pullback of abelian groups, which proves the equivalence of (i) and (ii).
The equivalence of (ii) and (iii) is explained by the following lemma.
\end{beweis}

\begin{lemma} \label{Postnikov-Faserungen}
Let $K\to L$ be a fibration between fibrant simplicial sets and $n\bth 0$. This map induces isomorphisms on homotopy groups in degrees $s\uber n$ for all basepoints if and only if for $s\bth n+2$ the induced maps
    $$ K_s \to M_sK\times_{M_sL}L_s $$
and the map
    $$ M_{n+2}K\to M_{n+2}Y\times_{Y_0}X_0$$
are surjective. Here the last map is induced by some map $\ast\to\partial\Delta^{n+2}$.
\end{lemma}

\begin{beweis}
This will be proved in \cite{Biedermann:n-types}.
\end{beweis}

Now we are going to determine the $n$-\mc{G}-cofibrant objects. Remember that cofibrant objects in the \mc{G}-structure coincide with the Reedy cofibrant ones.
\begin{lemma} \label{n-G-kofibrante Objekte}
\begin{punkt}
      \item An object $A^\bullet$ in $c\mc{M}$ is $n$-\mc{G}-cofibrant if and only if it is Reedy cofibrant and $\emph{\naturalpi{s}{A^\bullet}{G}}=0$ for all $G\aus\mc{G}$ and $s\uber n$. 
      \item An $n$-\mc{G}-cofibrant approximation functor is given by $Q_n=\sk_{n+1}\wt{\frei}$.
\end{punkt}
\end{lemma}

\begin{beweis}
Obvious from lemma \ref{n-G-Kofaserungen}.
\end{beweis}

\begin{Bem} \label{Kofas. zwischen n-G-kof. Objekten}
On $n$-\mc{G}-cofibrant objects the $n$-\mc{G}-structure and the \mc{G}-structure coincide.
\end{Bem}

\begin{Def} \label{Ko-Postnikov-Turm}
Let $X^\bullet$ be an object in $c\mc{M}$. The skeletal filtration of a Reedy cofibrant approximation to $X^\bullet$ consists of $n$-\mc{G}-cofibrant approximations $X^\bullet_n$ to $X^\bullet$ for the various $n$, and these assemble into a sequence
    $$ X^\bullet_0\to X^\bullet_1\to X^\bullet_2\to ... \to X^\bullet $$
which captures higher and higher natural homotopy groups. So this can be viewed as a {\bf Postnikov cotower} for $X^\bullet$. 
\end{Def}

\subsection{The tower of truncated homotopy categories}
Now we will study the homotopy categories associated to the truncated structures.

\begin{Bem} \label{volle Unterkategorie}
The functor $\id\co c\mc{M}^{\mc{G}}\to c\mc{M}^{n-\mc{G}}$ preserves weak equivalences and fibrations. It is therefore a right Quillen functor, whose left adjoint is given by $Q_n=\sk_{n+1}\wt{\frei}$. We have an induced pair of adjoint derived functors:
    $$ LQ_n\cong L(\id)\co\ho{c\mc{M}^{n-\mc{G}}}\leftrightarrows\ho{c\mc{M}^{\mc{G}}}:\!R(\id) $$
The unit $\id\to R(\id)L(\id)$ of this adjunction is a natural equivalence. Hence, 
    $$ LQ_n\cong L(\id)\co\ho{c\mc{M}^{n-\mc{G}}}\inj\ho{c\mc{M}^{\mc{G}}}$$
is an embedding of a full subcategory with a right adjoint given by $R(\id)$. In the same way, we can view $\id\co c\mc{M}^{(n+1)-\mc{G}}\to c\mc{M}^{n-\mc{G}}$ as a right Quillen functor. It induces a pair of adjoint derived functors
   $$ L(\id)\co\ho{c\mc{M}^{n-\mc{G}}}\leftrightarrows\ho{c\mc{M}^{(n+1)-\mc{G}}}:\!R(\id)=\sigma_n, $$
where the left adjoint is again an embedding of a full subcategory.
We obtain a tower 
    $$ ... \to \ho{c\mc{M}^{(n+1)-\mc{G}}} \toh{\sigma_n} \ho{c\mc{M}^{n-\mc{G}}}\to ... \to \ho{c\mc{M}^{1-\mc{G}}}\toh{\sigma_0} \ho{c\mc{M}^{0-\mc{G}}} $$
of categories, which can be identified with a tower of full subcategories of $\ho{c\mc{M}^{\mc{G}}}$ given by coreflections.
We can characterize the objects in $\ho{c\mc{M}^{n-\mc{G}}}$ viewed as a subcategory of $\ho{c\mc{M}^{\mc{G}}}$ by their natural homotopy groups. An object $X^\bullet$ is in the image of $\ho{c\mc{M}^{n-\mc{G}}}$ if and only if it is \mc{G}-equivalent to its $n$-\mc{G}-cofibrant replacement, i.e. if we have:  
\begin{equation*}
    \naturalpi{s}{X^\bullet}{G}=0 \hbox{ for } s\uber n
\end{equation*}
We have to relate all this to $\ho{\mc{M}}$ by the following statement, whose analogue for the Reedy structure is well known. The lemma is cited from \cite[Prop. 8.1.]{Bou:cos}.
\end{Bem}

\begin{lemma} \label{Tot-Delta-Quillen-Paar}
The functors 
    $$\xymatrix{ \mc{M} \ar@<2pt>[rr]^-{\frei\otimes^{\rm pro}\Delta^\bullet} & & c\mc{M}^{\mc{G}} \ar@<2pt>[ll]^-{\Tot} } .$$
form a Quillen pair.
\end{lemma}

\begin{Bem} \label{Dasselbe Quillenpaar}
The natural transformation $\frei\otimes^{\rm pro}\Delta^\bullet\to r^0$ gives a Reedy cofibrant replacement by \cite[16.1.4.]{Hir:loc} and hence a \mc{G}-cofibrant replacement. It follows that both induce the same left derived functor:
    $$\xymatrix{ \ho{\mc{M}} \ar@<2pt>[rrr]^-{Lr^0=\,\frei\otimes^{\rm pro}\Delta^\bullet} &&& \ho{c\mc{M}^{\mc{G}}} \ar@<2pt>[lll]^-{R\Tot} } .$$
We can look at the composition
    $$ \mc{M}\toh{r^0} c\mc{M}^{\mc{G}} \toh{\id} c\mc{M}^{n-\mc{G}}  $$
and the compostion of induced derived functors:
\begin{equation}\label{theta_n} 
      \xymatrix{\ho{\mc{M}}\ar[r]^-{Lr^0} & \ho{c\mc{M}^{\mc{G}}} \ar[r]^-{R(\id)} & \ho{c\mc{M}^{n-\mc{G}}} }.  
\end{equation}
\end{Bem}

\begin{Def} 
We will denote the composition \Ref{theta_n} of functors by \boldmath$\theta_n$\unboldmath. We arrive at the following diagram:
     $$\xymatrix@=20pt{ && \mc{T}=\ho{\mc{M}} \ar[d]^{\theta_n}\ar[dl]_{\theta_{n+1}}\ar[drr]^{\theta_0} && \\
        ... \ar[r] & \ho{c\mc{M}^{(n+1)-\mc{G}}} \ar[r]_-{\sigma_n} & \ho{c\mc{M}^{n-\mc{G}}}\ar[r] &  ... \ar[r]_-{\sigma_0} & \ho{c\mc{M}^{0-\mc{G}}}  } $$
This diagram is a $2$-commuting diagram of functors. For details on $2$-commutativity we refer to \cite{Hov:model}.
$2$-commutativity is provided by the relation $QQ\simeq Q$.
We call this {\bf the tower of truncated homotopy categories} associated to \mc{M} and \mc{G}. 
\end{Def}

\begin{Bem} 
Note that $\sk_{n+1}$ and $\id$ taken as functors $c\mc{M}^{\mc{G}}\to c\mc{M}^{n-\mc{G}}$ preserve weak equivalences and possess right adjoints, but they are not left Quillen functors since they do not preserve cofibrations, although $\sk_{n+1}$ and id preserve Reedy cofibrations. In fact $\sk_{n+1}$ does not map \mc{G}-cofibrations to \mc{G}-cofibrations. Therefore the functors $R(\id)\co\ho{c\mc{M}^{\mc{G}}}\to\ho{c\mc{M}^{n-\mc{G}}}$ and $\sigma_n\co\ho{c\mc{M}^{(n+1)-\mc{G}}}\to\ho{c\mc{M}^{n-\mc{G}}}$ do not have right adjoints.
\end{Bem}

\subsection{The \mc{G}-homotopy spectral sequence}

In \cite[2.9.]{Bou:cos} and \cite[5.8.]{Bou:cos} Bousfield describes, how a \mc{G}-resolution model structures leads to a \mc{G}-homotopy spectral sequence, whose convergence behaviour can be studied by a \mc{G}-completion functor. The spectral sequence $E_*^{*,*}\hrl X,Y\hrr$ in question is the Bousfield-Kan spectral sequence associated to the cosimplicial space $\map\hrl X,Y^\bullet\hrr$, where $X$ is cofibrant and $Y^\bullet$ is a \mc{G}-fibrant approximation to $Y$. Its differentials go like
   $$ d_r\co E_r^{s,t}\hrl X,Y\hrr\to E_r^{s+r,t+r-1}\hrl X,Y\hrr $$  
with $E_1^{s,t}\hrl X,Y\hrr\cong N^s\pi_t\map\hrl X,Y^\bullet\hrr$ and $E_2^{s,t}\hrl X,Y\hrr\cong\pi^s\pi_t\map\hrl X,Y^\bullet\hrr$, where $N^s$ denotes normalization. If we are in an unstable situation this spectral sequence is only defined partially for $0\sth s\sth t$ and is fringed along the line $s=t$. A thorough study of these phenomena is given in \cite{Bou:fringe}.

On the other hand we get an exact couple by applying $\stabhom{\frei}{Y^\bullet}_\mc{G}$ to our Postnikov cotower in \ref{Ko-Postnikov-Turm}. It follows from the dual version of \cite[3.9.]{GoHop:moduli}, that this is the derived exact couple of the above one. So up to reindexing we get from the $E_2$-term onwards an isomorphic spectral sequence. This is an alternative description of the \mc{G}-homotopy spectral sequence.

\end{document}